\documentclass[reqno,12pt]{article} \usepackage{fontenc} 
\usepackage{amsfonts}  \usepackage{amssymb}

\usepackage{pstricks} 

\def\alert#1{\smallskip{\hskip\parindent\vrule%
\vbox{\advance\hsize-2\parindent\hrule\smallskip\parindent.4\parindent%
\narrower\noindent#1\smallskip\hrule}\vrule\hfill}\smallskip}

\newcommand{\charr}{\mbox{char~}}
\newcommand{\SSS}{\mathcal{S}}

\newcommand{\RR}{\mathbb{R}}

\newcommand{\riff}{\Longrightarrow}

        \newcommand{\ZZ}{\mathbb{Z}}
\newcommand{\CC}{\mathbb{C}}

\newtheorem{thm}{Theorem}[section]
\newtheorem{lem}[thm]{Lemma}

\newtheorem{rem}[thm]{Remark}

\usepackage{amssymb}
\usepackage{amsfonts}
\def\alert#1{\smallskip{\hskip\parindent\vrule%
\vbox{\advance\hsize-2\parindent\hrule\smallskip\parindent.4\parindent%
\narrower\noindent#1\smallskip\hrule}\vrule\hfill}\smallskip}

\usepackage{pstricks}

\begin{document}



\begin{center}
{\bf \Large 
Irreducibility of the Cayley-Menger determinant,}\\
\vspace{.2cm}
 {\bf \Large and of a class of related polynomials}\\~~\\
{Mowaffaq Hajja$^{1}$, Mostafa Hayajneh$^2$, Bach  Nguyen$^3$, and Shadi Shaqaqha$^4$}\\~\\
(1): P. O. Box 1, Philadelphia University, 19392, Amman, Jordan\\
mowhajja@yahoo.com\\~~\\
(2), (4): Yarmouk University, Irbid, Jordan\\
hayaj86@yahoo.com,~shadish2@yahoo.com\\~~\\
(3): Louisiana State University, Baton Rouge, LA, USA\\
bnguy38@lsu.edu
\end{center}

\vspace{.4cm} \begin{quote}
{\small \bf Abstract.} {\small If $S$ is a given regular $n$-simplex, $n \ge 2$, of edge length $a$, then the distances $a_1$, $\cdots$, $a_{n+1}$ of an arbitrary point in its affine hull to its vertices are related by the fairly known elegant relation $\phi_{n+1} (a,a_1,\cdots,a_{n+1})=0$, where
$$\phi = \phi_t (x, x_1,\cdots,x_{n+1}) = \left( x^2+x_1^2+\cdots+x_{n+1}^2\right)^2
- t\left( x^4+x_1^4+\cdots+x_{n+1}^4\right).$$
The natural question whether this is essentially the only relation is answered positively by M. Hajja, M. Hayajneh, B.  Nguyen, and Sh. Shaqaqha 
in a recently submitted  paper entitled  {\it Distances from the vertices of a regular simplex}.
In that paper, the authors made use of the irreducibility of the polynomial $\phi $ in the case when $n \ge 2$, $t=n+1$, $x= a \ne 0$, and $k = \RR$,  but supplied no proof, promising to do so in another paper that is turning out to be this one. It is thus the main aim of this paper  to establish that irreducibility. In fact, we treat the irreducibility of $\phi$ without  restrictions on $t$, $x$, $a$, and $k$. As a by-product, we obtain new proofs of 
results pertaining to the irreducibility of the general Cayley-Menger determinant that are more general than those established by C. D'Andrea and M.  Sombra in {\it 
Sib. J. Math. {\bf 46}, 71--76.}
}
\end{quote}

\vspace{.4cm} \noindent {\bf Keywords:} Cayley-Menger determinant; circumscriptible simplex;   discriminant; homogeneous polynomial;  irreducible polynomial;      isodynamic simplex; isogonic simplex; orthocentric simplex; Pompeiu's theorem; pre-kites; quadratic polynomial; regular simplex;  symmetric polynomial; tetra-isogonic simplex;     volume of a simplex

\vspace{.4cm} \noindent 
\section{Introduction} \label{111}
Let $S=[A_1, \cdots, A_{n+1}]$, $n \ge 2$,  be a regular $n$-simplex of edge length $a$, and let $B$ be an arbitrary point in its affine hull. It is fairly well known that the distances $a_1,\cdots,a_{n+1}$  from 
$B$ to the vertices of $S$ satisfy the elegant relation
\begin{eqnarray} \label{f-a}
\left(a^2+a_1^2+\cdots+a_{n+1}^2\right)^2 &=& (n+1)\left(a^4+a_1^4+\cdots+a_{n+1}^4\right);
\end{eqnarray}
see \cite{Bentin-1}. The natural question whether this is (essentially) the only relation was posed and answered in the affirmative in \cite{Shadi}.
However, the proof rests heavily on a result whose proof was not included therein, but was postponed for another paper, namely this one. This result states that  the polynomial obtained from (\ref{f-a}) by thinking of $a_1,\cdots,a_{n+1}$ as indeterminates
is irreducible over  $\RR$  when $a \ne 0$ and  $n \ge 2$.

The main goal of this paper is to prove this, and to actually  give a complete treatment of the irreducibility of the more general polynomial 
\begin{eqnarray} \label{g}
g &=& \left(a^2+x_1^2+\cdots+x_{n+1}^2\right)^2 - t\left(a^4+x_1^4+\cdots+x_{n+1}^4\right) \in k[x_1,\cdots,x_{n+1}]
\end{eqnarray}
for any field $k$, any  $t \in k$,  and any $a \in k$ (including the case $a=0$).  It turns out that if $\charr k \ne 2$, then $g$ is irreducible except in the two cases  $(n,a,t)=(2,0,2)$ and $(n,a,t)=(2,0,3)$. When $(n,a,t)=(2,0,2)$, $g$ is the Heron polynomial with the well known factorization given in (\ref{uuu}). When $(n,a,t)=(2,0,3)$,
$g$ is  reducible if and only if $k$ contains a primitive third root of unity, in which case $g$ factors as in (\ref{vvvv}).
If $\charr k = 2$, then $g$ reduces to the polynomial $(1-t)(a+x_1+\cdots+x_{n+1})^4.$

The cases when $n=0$ and $n=1$ are ignored because they have no geometric significance. Also,  the case $n=1$ turned out to be rather lengthy and complex, but it gives rise to interesting  number theoretic aspects. We do not include these cases in this paper.

The polynomial obtained from $g$ by replacing $n+1$ with $n$ and $t$ with $n$ is essentially the Cayley-Menger determinant of a special  $n$-simplex, named a prekite, that was introduced and studied  in \cite{prekites}, and calculated in Theorem 4.1 there. If we put $a=0$ in $g$, then the instances $t=0$, $t=n$, $t=n-1$, and $t=n-2$ are very closely related to  the Cayley-Menger determinants, with respect to certain parametrizations,  of orthocentric, isodynamic, circumscriptible, and tetra-isogonic $n$-simplices; see Theorem 5.2 of \cite{impurity}. Thus the results of this paper would help give complete factorizations of 
the   Cayley-Menger determinants of the aforementioned five special families of $n$-simplices.
This raises the question whether the Cayley-Menger determinant $M$ of a general $n$-simplex is irreducible. It is proved in Theorem  \ref{CM} that  the answer is affirmative except in the single case $n=2$, in which case $M$ is the Heron polynomial given in (\ref{uuu}). This result has been obtained earlier in \cite{Siberian} when $k$ is $\RR$ or $\CC$. We should note, however, that the result in \cite{Siberian} does not imply  the results in this paper, and does not give any information about the irreducibility of the Cayley-Menger determinant of any of the five special families mentioned earlier. This is because the Cayley-Menger determinants of these special families 
are obtained from the general Cayley-Menger determinant $M$ given in (\ref{CMCM})
by replacing each  $x_{ij}^2$, $1 \le i < j \le n+1$, by quantities like $x_i+x_j$,  $x_i x_j$, $(x_i+x_j)^2$, $x_i^2+x_ix_j+x_j^2$; see \cite[Theorem ~5.1, p.~54]{impurity}. 
Such substitutions obviously do not preserve irreducibility. Furthermore, our treatment is not confined to the cases when $k = \RR$ or $\CC$.

\bigskip
The paper is organized as follows. Section \ref{222} contains  preliminary facts that we shall need, and freely use, in the sequel. These are quite few, and come from the  elementary theories of symmetric and of homogeneous polynomials. Section \ref{333} contains  proofs of the main results, namely Theorems \ref{DDD}, \ref{CCDD}, and \ref{CM}. These proofs make use of Lemmas \ref{AAA}, \ref{BBB}, and \ref{CCC}.

\section{Preliminaries} \label{222} In this section, we put together definitions and simple  facts pertaining to  symmetric and homogeneous polynomials that we shall freely use throughout the paper. For ease of reference, we also include a simple theorem about the factorization of quadratics.

Let $A$ be  any commutative ring with identity, and let $B=A[x_1,\cdots,x_n]$ be the polynomial ring in the $n$ indeterminates $x_1,\cdots,x_n$. Then the symmetric group  $\SSS_n$  acts  as permutations on the indices of $x_1,\cdots,x_n$, and hence as $A$-automorphisms on $B$. A polynomial in $B$ is called {\it symmetric} if it is invariant under the action of (every element of) $\SSS_n$.  It is well known that the ring $T=A[x_1,\cdots,x_n]^{\SSS_n}$ of symmetric polynomials is given by
\begin{eqnarray}\label{symsym}
T=A[x_1,\cdots,x_n]^{\SSS_n} &=&A[e_1,\cdots,e_n],
\end{eqnarray}
 where $e_1, \cdots, e_n$ are the so-called {\it elementary symmetric polynomials} defined by
$$e_j = \sum_{1 \le i_1 < i_2 < \cdots < i_j \le n}  x_{i_1} x_{i_2} \cdots x_{i_j};$$
see \cite[Theorem 6.1, p.~191]{Lang}. 

A non-zero polynomial $\Phi$ in $B$ is said to be {\it homogeneous} if all of its terms are of the same (total) degree $d$. More precisely, if $$\Phi (\lambda x_1, \cdots, \lambda x_n) = \lambda^d \Phi  
(x_1, \cdots, x_n)$$
for any variable $\lambda$.

It is clear that every non-zero polynomial $\Phi  \in B$ of (total) degree $n \ge 0$ can be written uniquely in the form $$\Phi  = h^{(0)} + h^{(1)} + \cdots + h^{(n)},$$
where  $h^{(j)}$ is either 0 or homogeneous of  degree $j$, and $h^{(n)}$ is not zero. The last term $h^{(n)}$ is called the {\it leading homogeneous component} of $h$, and will be denoted by $\Phi^*$. Notice that 
\begin{eqnarray}
\deg \Phi^* = \deg \Phi.
\label{degh}
\end{eqnarray}

If $A$ is an integral domain, then it is easy to see that 
\begin{eqnarray}
(\Phi \Psi)^* = \Phi^* \Psi^*.
\label{fgh}
\end{eqnarray}
Then it follows from (\ref{degh}) and (\ref{fgh}) that  if $f$ is reducible, then $f^*$ is reducible. It also follows that 
\begin{eqnarray}\label{hhh}
\mbox{every factor of a homogeneous polynomial is  homogeneous};
\end{eqnarray}
 see \cite[Theorem 10.5, p.~28]{Walker}. 
To see that this is not true if $A$ was not an integral domain,  consider the factorization $x^2 = (x-2)(x+2)$ in $\ZZ_4[x]$.

We end this section by proving a  simple property of quadratic polynomials that we shall use later.

\begin{lem} \label{Q}
Let $A$ be an integral domain with $\charr A \ne 2$, and let $A[x]$ be the polynomial ring over $A$ in the indeterminate $x$. Let $H =ax^2+bx+c \in A[x]$, with $a \ne 0$, and let 
$\Delta = b^2 - 4 ac$ be its discriminant. Then
\begin{itemize}
\item[(i)] If $H$ is reducible, then $\Delta$ is a square (in $A$).
\item[(ii)] If $\Delta$ is a square, and if $2a$ is a unit, then $H$ is reducible.
 \item[(iii)] If $H$ is a square, then $\Delta = 0$.
 \item[(iv)] If $\Delta = 0$, and $a$ is a square and $2a$ is a unit, then  $H$ is a square.
\end{itemize}
\end{lem}

\vspace{.2cm} \noindent {\it Proof.} To prove (i), suppose that $H$ is reducible. Then $H=(U x + V)(ux+v)$ for some $U, V, u, v \in A$. Multiplying out and equating coefficients, we obtain $Uu=a$, $Uv+Vu=b$, and $Vv=c$.
Therefore 
$\Delta = b^2 - 4ac =  (Uv+Vu)^2 - 4 (Uu)(Vv) =(Uv-Vu)^2$,
which is a square, as desired.

To prove (ii), suppose that $\Delta = b^2 - 4ac$ is a square, say $\Delta = \delta^2$, and that $a$ is a unit. Letting
$\alpha = \left(-b+\delta\right)\left(2a\right)^{-1}$ and $\beta = \left(-b-\delta\right)\left(2a\right)^{-1}$,
we can easily see that $H=a(x-\alpha)(x-\beta),$ and hence is reducible, as desired. 

To prove (iii), suppose that $H$ is a square, say $H=(ux+v)^2$, where $u, v \in A$. Then $a = u^2$, $b=2uv$, and $c=v^2$.
Thus  $b^2-4ac=0$.

To prove (iv), suppose that  $2a$ is a unit and that $a$ is a square, say $a=u^2$, and $b^2-4ac=0$. Then it is easy to check that 
$H=\left(ux+ b (2u)^{-1} \right)^2,$
a square.  \hfill $\Box$

\begin{rem} 
{\em
 The assumption  that $2a$ is a unit in (ii) and (iv) above
cannot be replaced by the weaker assumption that $a$ is a unit. For (ii), consider the example $A = \ZZ[\sqrt{5}]$ and  $H=x^2+\sqrt{5} x + 1$. Then $a=\Delta = 1$, but $H$ is irreducible. In fact, if $H=(x+s)(x+t)$, then $st=1$, $s+t=\sqrt{5}$, and $(s-t)^2 = (s+t)^2-4st=1,$
and hence $s-t = \pm 1$. In view of the fact that  $s+t=\sqrt{5}$, we obtain
$2s=\sqrt{5}\pm 1$, contradicting the fact that $s \in A$.
For  (iv), let $D = \{ f(x) = c_0+c_1 x+\cdots +c_n x^n\in \ZZ[x] : c_1 \mbox{~is even}\},$
and let $H = x^2 + 2x + 1$. Then $\Delta = 0$, and $H$ is not a square in $D$.}
\end{rem}

\section{The main results} \label{333}
In this section, we establish, in Theorems \ref{CCDD}, \ref{DDD}, and \ref{CM}, the  main irreducibility theorems of this paper.  Lemmas \ref{AAA},  \ref{BBB},  and \ref{CCC} are  needed in the proofs. Lemma \ref{CBA} disposes of the case when $t=0$ in (\ref{g}).

\begin{lem} \label{CBA}
Let $k$ be a field,  and let  $g=c_0 + c_1x_1^2+\cdots + c_n x_n^2
\in k[x_1,\cdots, x_n]$, where $n \ge 1$ and $c_1,\cdots,c_n \in k \setminus \{0\}$. For $0 \le i \le n$, let  $t_i = c_i/c_n$, and let 
$h=t_0 + t_1x_1^2+\cdots + t_n x_n^2.$
Then $h$ (and hence  $g$) is reducible if and only if 
$-t_j$ is a square for all $j$, $0 \le j \le n-1$ and if
\begin{itemize}
\item[(1)] either $\charr k =2$, in which case
	 	\begin{eqnarray*}
	h&=&\left(\sqrt{t_0} + \sqrt{t_1} x_1+ \cdots +
\sqrt{t_{n-1}} x_{n-1} +  x_{n}          \right)^2,
\end{eqnarray*}
\item[(2)] or $\charr k \ne 2$ and  $n=1$, in which  case,
\begin{eqnarray*}
	h&=&\left(x_1 - \sqrt{-t_0}\right)\left(x_1 + \sqrt{-t_0}\right),
	\end{eqnarray*}
\item[(3)] or $\charr k \ne 2$,  $n=2$, and $t_0 = 0$, in which case 
\begin{eqnarray*}
	h&=&\left(x_2 - \sqrt{-t_1}x_1 \right)\left(x_2 + \sqrt{-t_1} x_1\right).
	\end{eqnarray*}
\end{itemize}
\end{lem}

\noindent {\it Proof.} 
Suppose that $h$ is reducible. Then any factorization of $h$ must be of the form
\begin{eqnarray} \label{Must-1}
h &=&(a_0 + a_1 x_1+ \cdots+a_nx_n)(b_0 + b_1 x_1+ \cdots+b_nx_n).
\end{eqnarray}
 We may also assume that $a_n=b_n=1,$ since $a_n b_n = 1$.
Comparing coefficients of $x_nx_i$, $1 \le i \le n-1$, and coefficients of $x_n$, we see that $b_i=-a_i$
for $0 \le i \le n-1$.
Therefore (\ref{Must-1}) can be rewritten as
\begin{eqnarray} \label{Must-2}
h &=&(a_0 + a_1 x_1+ \cdots+a_{n-1} x_{n-1}+ x_n)(-a_0 -a_1 x_1- \cdots -a_{n-1} x_{n-1} + x_n).
\end{eqnarray}
Comparing coefficients of $x_i^2$, $1 \le i \le n-1$, and the constant terms, we see that $-t_i = a_i^2$ for $0 \le i \le n-1$.
   Thus $- t_i$ is a square for all $i$, $0 \le i \le n-1.$ Of course, $t_n = 1$.
   
   If char ($k$) = 2, then $h$ factors into
\begin{eqnarray*}
h &=&(a_0 + a_1 x_1+ \cdots+a_{n-1} x_{n-1}+ x_n)^2.
\end{eqnarray*}

If char ($k$) $\ne 2$, and $n > 2$, then comparing the coefficients of $x_1x_2$ (in (\ref{Must-2})), we obtain $-2 a_1 a_2 = 0$, contradicting the assumptions that  $a_1^2 = -t_1 \ne 0$,  $a_2^2 = -t_2 \ne 0$, and $2 \ne 0$.

If (char ($k$) $\ne 2$, and) $n = 1$, then 
$h$ factors into $h = (a_0 + x_1)(-a_0 + x_1)$.

If char ($k$) $\ne 2$, and $n = 2$, then comparing the coefficients of $x_1$, we obtain $-2 a_0 a_1 = 0$, and hence $a_0 = 0$, since
$a_1^2 = -t_1 \ne 0$ and $2 \ne 0$.
Therefore $h$ factors into
$h = (a_1 x_1 + x_2)(-a_1 x_1 + x_2).$ This completes the proof. \hfill $\Box$

\bigskip It follows, for example, that the polynomial $x_1^2 + \cdots + x_n^2$ is reducible over the field $\CC$ of complex numbers if and only if $n=1$ or $n=2$. In particular, $x^2+y^2+z^2$ is irreducible over $\CC$, a fact that appears as a problem in  \cite{AMATYC}.

\bigskip 
In view of the  lemma above, we may exclude the case $t=0$ in (\ref{g}). Also, it is obvious that if $\charr k = 2$, then $g = (1-t) (a+x_1+\cdots+x_n)^4$.  That is why  the lemmas and theorems below  assume  that $t\ne 0$ and  $\charr k \ne 2$.

\begin{lem} \label{AAA}
Let $k$ be a field with $\charr k \ne 2$, and let  
\begin{eqnarray*}f &=& (x_1^2+\cdots+x_n^2)^2-t(x_1^4+ \cdots +x_n^4) \in k[x_1,\cdots,x_n], \end{eqnarray*}
where $t \in k$ and $t \ne 0$, and where $n \ge 3$. Then
$f$ has a factor that is symmetric in two of the variables if and only if $n=3$ and $t=2$. In this case, 
\begin{eqnarray}f &=& (x_1^2+x_2^2+x_3^2)^2-2(x_1^4+ x_2^4+ x_3^4)\nonumber \\&=&
(x_1+x_2+x_3)(-x_1+x_2+x_3)(x_1-x_2+x_3)(x_1+x_2-x_3). \label{uuu}
 \end{eqnarray}
 \end{lem}

\vspace{.3cm} \noindent {\it Proof.} 
We shall freely use the fact that factors of a homogeneous polynomial are homogeneous; see (\ref{hhh}). 

Let 
	\vspace{-.4cm}
	\begin{eqnarray}
&&y = x_{n-1}, ~ z = x_n,~ u = y+z,~v=yz, \label{yzuv}\\
&&S_2=x_1^2+\cdots+x_{n-2}^2,~S_4=x_1^4+\cdots+x_{n-2}^4,\label{S2S4}\\
&&R=k[x_1,\cdots,x_n] = k[x_1,\cdots,x_{n-2},y,z],~
 R_0=k[x_1,\cdots,x_{n-2},u,v]. \label{RR}
\end{eqnarray}
Clearly, $R$ and $R_0$ are polynomial rings over $k$ in the $n$ respective variables, and $R_0$ is the subring of $R$ consisting of the $(y,z)$-symmetric polynomials, i.e., polynomials  that are symmetric in the variables $y$ and $z$; see (\ref{symsym}). 

Using the identities
\begin{eqnarray*}
y^2+z^2&=&u^2-2v,\\
y^4+z^4&=&(y^2+z^2)^2-2v^2~=~(u^2-2v)^2-2v^2,
\end{eqnarray*}
and the definition
$$f=(S_2+y^2+z^2)^2-t(S_4+y^4+z^4),$$
we obtain (after few lines of computations, or by Maple)
\begin{eqnarray} \label{fSuv}
f&=&(4-2t)v^2 -4 ((1-t)u^2+S_2)v +S_2^2 + (1-t) u^4  + 2S_2 u^2 -t S_4.
\end{eqnarray}
Letting $F \in R_0$  denote the right hand side of (\ref{fSuv}), i.e.,
\begin{eqnarray}
F&=&  (4-2t)v^2 -4 ((1-t)u^2+S_2)v +S_2^2 + (1-t) u^4  + 2S_2 u^2 -t S_4, \label{F}
\end{eqnarray}
we see that  
$f$ has a $(y,z)$-symmetric factor if and only if $F$ is reducible in $R_0$.
From now on, we will be working in $R_0$.

Suppose that $F$ is reducible.

If $t=2$, then $F$ simplifies into
\begin{eqnarray*}
F_0&=& -4 (-u^2+S_2)v +S_2^2 - u^4  + 2S_2 u^2 -2 S_4 \\
&=&  4 (u^2-S_2)v -  (u^2-S_2)^2    + 2S_2^2 -2 S_4.
\end{eqnarray*}
Since $F_0$ is linear in $v$, it follows that every  factor of $F_0$ is a factor of $(u^2-S_2)$ and of $-(u^2-S_2)^2+2S_2^2 -2 S_4$.
But any factor of $u^2-S_2$ must contain $u$, while $2S_2^2 -2 S_4$ does not. Therefore $F_0$ is reducible if and only if $2S_2^2 -2 S_4 = 0$, i.e., if and only if $n=3$, in which case
 $S_4=x_1^4=(x_1^2)^2=S_2^2$. In this case $F_0$ is as given in (\ref{uuu}).

If  $t \ne 2$, then $F$ is a quadratic in $v$ with coefficients in $k[x_1,\cdots,x_{n-2},u]$. Also, its leading coefficient $(4-2t)$ is a unit. By (i) and (ii) of Lemma \ref{Q}, it is reducible  if and only if its discriminant $\Delta$ in $v$ is a square in $k[x_1,\cdots,x_{n-2},u]$. Using Maple, we find that
\begin{eqnarray} \label{D}
\Delta &=&  8t((t-1)u^4-2S_2u^2+(2-t)S_4+S_2^2).\end{eqnarray}
If $t=1$, then $\Delta = 8(-2S_2 u^2 + S_4 + S_2^2)$. By (iii) of Lemma \ref{Q}, this  cannot be a square, since its discriminant in $u$ is $8^3S_2(S_4 + S_2^2) \ne 0$.
Thus suppose that 
$$t \ne 1.$$ Then for $\Delta$ to be a square, we must have 
\begin{eqnarray} \label{DD}
\Delta &=&  8t((t-1)u^4-2S_2u^2+(2-t)S_4+S_2^2) ~=~ (\alpha u^2 + \beta u + \gamma)^2, \end{eqnarray}
where $\alpha$, $\beta$, and $\gamma$ belong to $k[x_1,\cdots,x_{n-2}]$.
Equating the coefficients of $u^4$ and $u^3$ in (\ref{DD}), we see that
$\alpha^2 = 8t (t-1) \ne 0$ and $2 \alpha \beta = 0$. Therefore $\beta =0$, and we have
\begin{eqnarray} \label{DDDD} 
\Delta &=&  8t((t-1)u^4-2S_2u^2+(2-t)S_4+S_2^2) ~=~ (\alpha u^2 + \gamma)^2. 
\end{eqnarray}
Equating coefficients in (\ref{DDDD}), we find that
$$\alpha^2 = 8t(t-1),~~\gamma^2 = 8t ((2-t)S_4+S_2^2),~~2 \alpha \gamma = 8t
 (-2S_2).$$
 It follows that $(t-1) ((2-t)S_4+S_2^2) =  (-S_2)^2.$
 This simplifies into $(2-t)((t-1)S_4  - S_2^2) =  0.$ Since $t \ne 2$, it follows that $S_2^2$ and $S_4$ are linearly dependent over $k$. This happens if and only if $n=3$, in which case $S_4=x_1^4=(x_1^2)^2=S_2^2$, and $t-1=1$. This is the case treated earlier. Thus $\Delta$ is not a square, a contradiction. This completes the proof. \hfill $\Box$

\begin{lem} \label{BBB}
Let $k$ be a field with $\charr k \ne 2$, and let  
\begin{eqnarray*}f &=& (x_1^2+\cdots+x_n^2)^2-t(x_1^4+ \cdots +x_n^4) \in k[x_1,\cdots,x_n], \end{eqnarray*}
where $t \in k$ and $t \ne 0$, and where $n \ge 3$. Then
$f$ has a linear factor if and only if $n=3$ and $t=2$. 
In this case, $f$ factors as in (\ref{uuu}).
\end{lem}

\vspace{.2cm} \noindent {\em Proof.} If $(n,t)=(3,2)$, then $f$ factors as in (\ref{uuu}), and we are done. So we assume that $(n,t)\ne (3,2)$,  that 
$f$ has a linear factor, say $g=a_1x_1+\cdots+a_nx_n$, and we seek a contradiction.

Clearly we may assume that $a_1 \ne 0$. By Lemma \ref{AAA},  $g$ is not  symmetric in the variables $x_2$ and $x_3$. Therefore $a_2 \ne a_3$.
 If $g_1$ is the polynomial obtained from $g$ by interchanging $x_2$ and $x_3$, then $g_1$ is not an associate of $g$, and it is also a factor of $f$. Therefore $gg_1$ is a $(x_2,x_3)$-symmetric factor of $f$ (of degree 2). This contradicts Lemma \ref{AAA}. \hfill $\Box$

\begin{lem} \label{CCC}
Let $k$ be a field with $\charr k \ne 2$, and let  
\begin{eqnarray*}f &=& (x_1^2+x_2^2+x_3^2)^2-t(x_1^4+x_2^4+x_3^4) \in k[x_1,x_2,x_3], \end{eqnarray*}
where $t \in k$ and $t \ne 0$. If $t=2$, then $f$ factors as in (\ref{uuu}).
If $t \ne 2$, then
$f$ is reducible if and only if $t=3$ and $k$ contains a primitive third root $\omega$ of 1. In this case, 
\begin{eqnarray}
	f&=&(x_1^2+x_2^2+x_3^2)^2 - 3(x_1^4+x_2^4+x_3^4) 	\nonumber \\
	&=& (-2)(x_1^2+\omega x_2^2 + \omega^2 x_3^2)(x_1^2+\omega^2 x_2^2 + \omega x_3^2). \label{vvvv}
\end{eqnarray}
Also, the factors in the right hand sides of (\ref{vvvv})
are irreducible.
 \end{lem}

\vspace{.2cm} \noindent {\em Proof.} There is nothing to prove in the case when $t=2$, since this is covered in Lemma \ref{BBB}. Thus we assume that $t \ne 2$.

Let $x_1, x_2, x_3$ be renamed as $x, y, z$, and suppose that $f$ is reducible and that $g$ is an irreducible factor of $f$. Since $t \ne 2$, Lemma \ref{BBB} implies that $g$ is not linear. Thus $g$ is  an irreducible quadratic. 
Thus
\begin{eqnarray*}
g&=&ax^2+by^2+cz^2+\alpha yz + \beta zx + \gamma xy,
\end{eqnarray*}
where $a, b, c, \alpha, \beta, \gamma$ are in $k$.
Also,
\begin{eqnarray*}
f&=& (x^2+y^2+z^2)^2 - t(x^4+y^4+z^4)\\
&=& (1-t)(x^4+y^4+z^4)+2(x^2y^2+y^2z^2+z^2x^2).
\end{eqnarray*}
Letting $s$ be the permutation $s = (x \mapsto y \mapsto z \mapsto x)$, we see that $f$ is divisible by $g$, $s(g)$, and $s^2 (g)$. Since $\deg (g s(g) s^2(g)) = 6 > \deg f$, and since $g, s(g), s^2(g)$ are irreducible, it follows that two (and hence all) of the polynomials $g$, $s(g)$, and $s^2 (g)$ 
are associates (i.e., constant multiples of each other). Thus
$g = \lambda s(g)$ for some $\lambda \in k$. Since
\begin{eqnarray*}
s(g)&=&ay^2+bz^2+cx^2+\alpha zx + \beta xy + \gamma yz,
\end{eqnarray*}
  it follows that
   $$b = \lambda a,~ c = \lambda b,~a = \lambda c,~
   \beta = \lambda \alpha,~\gamma = \lambda \beta,~
   \alpha = \lambda \gamma.$$
   Thus
   $$a = \lambda^3 a,~\alpha = \lambda^3 \alpha.$$
      If both $a$ and $\alpha$ are zero, then  $a=b=c=\alpha = \beta =  \gamma=0$, and $g=0$, a contradiction.
      Thus either $a \ne 0$ or $\alpha \ne 0$. In both cases $
      \lambda^3 = 1$.
      If $\lambda = 1$, then $a=b=c$ and 
      $\alpha = \beta =  \gamma$, and hence $g$ is symmetric, contradicting Lemma \ref{AAA}. 
      Thus $\lambda \ne 1$ and $\lambda^3 = 1$. Therefore $\lambda$ is a primitive third root of 1 (and $\charr k$ cannot be 3). Thus $g$ is of the form
      $$g = a (x^2+\lambda y^2+ \lambda^2 z^2) + \alpha (yz+\lambda zx+ \lambda^2 xy),$$
where $\lambda$ is a primitive third root of 1.
Applying the permutation $y \mapsto z \mapsto y$ to $g$, we obtain another factor
      $$h = a (x^2+\lambda^2 y^2+ \lambda z^2) + \alpha (yz+\lambda^2 zx+ \lambda xy),$$
of $f$ that is not associate of $g$. Therefore $gh$ divides $f$, and has the same degree as $f$. Therefore $f = c gh$, where $c \in k$.
Equating the coefficients of $x^3z$ and $x^2yz$ in the identity $f=cgh$, we obtain $-a\alpha c = 0$ and $2a\alpha  c - \alpha^2 c  = 0$. Thus $\alpha =0$ and
      $$g = a (x^2+\lambda y^2+ \lambda^2 z^2),~~
      h = a (x^2+\lambda^2 y^2+ \lambda z^2).$$
Therefore
\begin{eqnarray} gh &=& a^2 [(x^4+y^4+z^4) - (x^2y^2+y^2z^2+z^2x^2)].
\end{eqnarray}
But
\begin{eqnarray} f &=& (1-t) (x^4+y^4+z^4) +2(x^2y^2+y^2z^2+z^2x^2).
\end{eqnarray}
Therefore it follows from $f = c gh$ that
$(1-t)=  ca^2$ and $-a^2c = 2$. Hence  $1-t=  -2$ and $t=3$. Therefore
	\begin{eqnarray*}
	f&=&(x^2+y^2+z^2)^2-3(x^4+y^4+z^4)  \\
	&=&(-2)(x^2+\lambda y^2 + \lambda^2 z^2)(x^2+\lambda^2 y^2 + \lambda z^2),
\end{eqnarray*}
as desired. 

The two factors in (\ref{vvvv}) are irreducible because 
$f$ has no linear factors, by Lemma \ref{BBB}. This completes the proof. \hfill $\Box$

\begin{thm} \label{CCDD}
Let $k$ be a field with $\charr k \ne 2$, and let $a$ be a non-zero element of $k$. Let  
\begin{eqnarray*}g &=& (a^2 + x_1^2+\cdots+x_n^2)^2-t(a^4 +x_1^4+ \cdots +x_n^4) \in k[x_1,\cdots,x_n], \end{eqnarray*}
where $t \in k$ and $t \ne 0$, and where $n \ge 3$. Then $g$ is irreducible.
\end{thm}

\vspace{.2cm} \noindent {\em Proof.} We start with the case $n=3$, and we rename $x_1$, $x_2$, $x_3$ as $x$, $y$, $z$. Thus we are to prove that the polynomial
$$g = (a^2 + x^2 + y^2 + z^2)^2 - t(a^4 + x^4 + y^4 + z^4),~~a\ne 0, ~t \ne 0, ~\charr k \ne 2,$$ is irreducible. The general case will follow easily as shown later.

Suppose that $g$ is reducible, and that $g = \alpha \beta$, where $\alpha$ and $\beta$ are non-constant polynomials in $k[x,y,z]$. Let $G$, $A$, and $B$ be the leading homogeneous components of $g$, $\alpha$, and $\beta$, respectively. Then 
$$G=(x^2+y^2+z^2)^2-t(x^4+y^4+z^4),$$
and $G = AB$, and $A$ and $B$ are non-constant. Therefore $G$ is reducible. By Lemma \ref{CCC}, we have the following two cases:

\vspace{.2cm} {\bf Case 1.}
\begin{eqnarray}
G&=&(x^2+y^2+z^2)^2-2(x^4+y^4+z^4)\nonumber \\
&=&(x+y+z)(-x+y+z)(x-y+z)(x+y-z), \label{GUV}
\end{eqnarray}
and hence 
\begin{eqnarray*}
g&=&(a^2+x^2+y^2+z^2)^2-2(a^4+x^4+y^4+z^4).
\end{eqnarray*}

Suppose that  one of the factors $\alpha$ or $\beta$, say $\alpha$, is linear. Then we may assume that  either
$\alpha=x+y+z+c$ or $\alpha=x+y-z+c$, where $c \in k$.
In the first case, we plug  $(x,y,z)=(x,-x,-c)$ in $g = \alpha \beta$, and in the second case, we plug
$(x,y,z)=(x,-x,c)$. In both cases, we obtain
\begin{eqnarray*}
0&=&(a^2+2x^2+c^2)^2-2(a^4+2x^4+c^4)
\\&=& 4x^2(a^2+c^2) +(-a^4-c^4+2a^2c^2)\\
&=& 4x^2(a^2+c^2)-(a^2-c^2)^2.
\end{eqnarray*}
Therefore $a^2+c^2=a^2-c^2=0$, and hence $a=c=0$, contradicting the assumption that $a \ne 0$.

If $\alpha$ and $\beta$ are  quadratic irreducible, then  $A$ and $B$ are quadratic. Also $G=AB$ by (\ref{fgh}).  By  (\ref{GUV}),
we may assume that $A$ or $B$, say $A$,  is $(x+y+z)(-x+y+z)$. Thus
\begin{eqnarray*}
\alpha&=&(x+y+z)(-x+y+z)+L,
\end{eqnarray*}
where $L$ is a  linear polynomial. By applying the permutation $\sigma : x \mapsto y \mapsto x$ and then $\tau : x \mapsto z \mapsto x$, we see that the two (irreducible) polynomials
\begin{eqnarray*}
\alpha_1=(x+y+z)(x-y+z)+\sigma (L) \mbox{~~and~~} \alpha_2=(x+y+z)(x+y-z)+ \tau (L)
\end{eqnarray*}
are also factors of $g$. Since  the coefficients of $x^2$ in the polynomials  $\alpha_1$ and $\alpha_2$ are the same (and equal 1), and since they are not equal, it follows that $\alpha_1$ and $\alpha_2$   are not associates.
By considering the coefficients of $y^2$ and $z^2$, we see that no two of 
the polynomials  $\alpha$, $\alpha_1$ and $\alpha_2$ are associates.
Therefore $\alpha \alpha_1 \alpha_2$ divides $g$, a contradiction since $\deg g = 4 < 6$.

Therefore $g$ cannot be reducible.

\vspace{.2cm}
{\bf Case 2.}
\begin{eqnarray*}
G&=&(x^2+y^2+z^2)^2-3(x^4+y^4+z^4)\\
&=&(-2)(x^2+\omega y^2+ \omega^2 z^2)(x^2+\omega^2 y^2+ \omega z^2),
\end{eqnarray*}
where $\omega \in k$ is a primitive third root of 1, and where the quadratics on the right hand side are irreducible.
Therefore $A$ and $B$ are the polynomials $$u(x^2+\omega y^2+\omega^2 z^2) \mbox{~and~} v(x^2+\omega^2 y^2+\omega z^2),$$
where $$u, v \in k,~uv=-2.$$
Hence $\alpha$ and $\beta$  are the polynomials $$u(x^2+\omega y^2+\omega^2 z^2 +L) \mbox{~and~} v(x^2+\omega^2 y^2+\omega z^2 +K),$$
where $L$ and $K$ are linear polynomials. Also, 
\begin{eqnarray}
g&=&(a^2+x^2+y^2+z^2)^2-3(a^4+x^4+y^4+z^4) \nonumber\\
&=&(-2)(x^2+\omega y^2+ \omega^2 z^2+L)(x^2+\omega^2 y^2+ \omega z^2+K). \label{Kh1}
\end{eqnarray}
Let $L_0$ and $K_0$ be the constant terms of  $L$ and $K$, respectively. Plugging $x=y=z=0$ in  (\ref{Kh1}), we obtain
$-2 L_0 K_0 =  -2 a^4$, and hence 
\begin{eqnarray}
 L_0 K_0 &\ne& 0. \label{Kh2}
\end{eqnarray}

Applying the permutation $\sigma : x \mapsto y \mapsto x$, we obtain
\begin{eqnarray*}
g&=&(a^2+x^2+y^2+z^2)^2-3(a^4+x^4+y^4+z^4)\\
&=&(-2)(y^2+\omega x^2+ \omega^2 z^2+ \sigma (L))(y^2+\omega^2 x^2+ \omega z^2+\sigma (K)),\\
&=&(-2)(x^2+\omega^2 y^2+ \omega z^2+ \omega^2 \sigma (L))(x^2+\omega y^2+ \omega^2 z^2+\omega \sigma (K)).
\end{eqnarray*}
Since $k[x,y,z]$ is a unique factorization domain, it follows that 
\begin{eqnarray}
\sigma (L) = \omega K,~ \sigma (K) = \omega^2 L. \label{Kh3}
\end{eqnarray}
Similarly, if $\tau$ is the transposition $x \mapsto z \mapsto x$, then
\begin{eqnarray}
\tau (L) = \omega^2 K,~ \tau (K) = \omega L. \label{Kh4}
\end{eqnarray}
Observing that the constant terms of $L$ and $K$ are unchanged under permutations on $x, y, z$, and using  
(\ref{Kh3}) and (\ref{Kh4}), we obtain
\begin{eqnarray}
L_0 = \omega K_0,~L_0 = \omega^2 K_0.
\end{eqnarray}
It follows that $L_0 = K_0 =0$, contradicting (\ref{Kh2}).

Thus we have proved the case $n=3$. If $n \ge 3$, we let $h$ be obtained from $g$ by putting $x_j = 0$ for all $j \ge 4$. Since $h$ is irreducible by the case $n=3$, it follows that $g$ is irreducible. In fact, if $g=\alpha \beta$, and if $A$ and $B$ are obtained from $\alpha$ and $\beta$ by plugging $x_j = 0$ for all $j \ge 4$, then either $A$ is zero or $A$ is homogeneous of the same degree as $\alpha$. Since $h=AB$, and $h \ne 0$, it follows that $A\ne 0$ and therefore $h$ is reducible, a contradiction. 

This completes the proof.
\hfill $\Box$

\begin{thm} \label{DDD}
Let $k$ be a field with $\charr k \ne 2$, and let  
\begin{eqnarray*}f &=& (x_1^2+\cdots+x_n^2)^2-t(x_1^4+ \cdots +x_n^4) \in k[x_1,\cdots,x_n], \end{eqnarray*}
where $t \in k$ and $t \ne 0$, and where $n \ge 3$. Then
$f$ is reducible in and only in the following  two cases:
\begin{itemize}
	\item[\rm{(i)}] $n=3$ and $t=2$, in which case $f$ is as given in (\ref{uuu}). 
			\item[\rm{(ii)}] $n=3$,  $t=3$, and $k$ contains a primitive third root $\omega$ of unity, in which case $f$ is as given in (\ref{vvvv}).
	\end{itemize}
\end{thm}

\vspace{.2cm} \noindent {\em Proof.} The  case when $n=3$ was completely treated in  Lemmas \ref{BBB} and \ref{CCC}. So we assume that $n \ge 4$.
For any polynomial $F \in k[x_1,\cdots,x_n]$, let $F^* \in k[x_1,x_2,x_3]$ be the polynomial obtained from $F$ by putting $x_4=1$ and $x_j = 0$ for all $j \ge 5$. Then $f^*$ is the case $n=3$ and $a=1$ of Theorem \ref{CCDD}, and is hence irreducible. To show that  $f$ is irreducible, we suppose that  $f = gh$, and we show that $g$ or $h$ is a constant.

Since $f^*=g^*h^*$, and since $f^*$ is irreducible, it follows that  $g^*$ or $h^*$, say $g^*$, is a constant.
Since $f^* \ne 0$, it follows that $g^*$ is a non-zero constant. Therefore 
$$4 = \deg (f^*) = \deg (g^*) + \deg (h^*)= 0 + \deg (h^*).$$
Hence $\deg (h^*) = 4$. Since $\deg (h) \ge \deg (h^*)$, it follows that $\deg (h) = 4$. Therefore $4 = \deg (f) = \deg (g) + \deg (h) = \deg (g) + 4$, and hence $g$ is a constant.

Thus $f$ is irreducible,  as desired. \hfill $\Box$

\begin{thm} \label{CM}
Let $k$ be a field with $\charr k \ne 2$, and let 
$R=k[x_{ij} : 1 \le i < j \le n+1]$, $n \ge 2$, be the polynomial ring over $k$ in the set $X=\{x_{ij} : 1 \le i < j \le n+1\}$ of $(n+1)n/2$ indeterminates.
For $1 \le i \le j \le n+1$, let us make the convention that $x_{i,j} = x_{j,i}$ and $x_{j,j} = 0$.  Let $M$ be the Cayley-Menger determinant in $X$, i.e.,
 $M$ is the $(n+2)\times (n+2)$ determinant whose entries $c_{i,j}$, $0 \le i, j \le n+1$, are given by
 \begin{eqnarray*}
 c_{i,j} &=& \left\{ \begin{array}{lcl} 0 &\mbox{if}& i=j,\\  
                                 1 &\mbox{if}& i \ne j  \mbox{~and~} ij=0,\\
 x_{i,j}^2 & \mbox{otherwise.} \end{array}\right.
 \end{eqnarray*}
Thus
\begin{eqnarray} \label{CMCM}
M&=& \left|
\begin{array}{cccccc}
0&1&1&1& \cdots\cdots & 1\\
1&0&x_{1,2}^2&x_{1,3}^2& \cdots\cdots&x_{1,n+1}^2\\
1&x_{2,1}^2&0&x_{2,3}^2& \cdots\cdots&x_{2,n+1}^2\\
1&x_{3,1}^2&x_{3,2}^2& 0& \cdots\cdots&x_{3,n+1}^2\\
\cdots \cdots & \cdots \cdots & \cdots \cdots & \cdots \cdots & \cdots \cdots & \cdots \cdots \\
\cdots \cdots & \cdots \cdots & \cdots \cdots & \cdots \cdots & \cdots \cdots & \cdots \cdots \\
\cdots \cdots & \cdots \cdots & \cdots \cdots & \cdots \cdots & \cdots \cdots & \cdots \cdots \\
1&x_{n+1,1}^2&x_{n+1,2}^2& x_{n+1,3}^2& 
\cdots \cdots & 0
\end{array}
\right|.
\end{eqnarray}
Then 
\begin{enumerate}
	\item[(i)] $M$ is homogeneous of homogeneity degree $2n$,
	\item[(ii)] $M$ is reducible if and only if $n=2$. In this case, if we let 
$$x_{1,2} = z,~x_{2,3} = x,~x_{1,3}=y,$$
then
\begin{eqnarray}
M&=& \left|
\begin{array}{cccccc}
0&1&1&1\\
1&0&z^2&y^2\\
1&z^2&0&x^2\\
1&y^2&x^2& 0
\end{array}
\right| \nonumber \\ &=& -(x + y + z) (-x+y+z) (x - y + z) (x+y-z). \label{Heron-M}
\end{eqnarray}
\end{enumerate}
\end{thm}

\vspace{.2cm} \noindent {\em Proof.} 
(i) Let $M_1$ be the matrix obtained from $M$ by replacing each $x_{i,j}$ in $M$ by $\lambda x_{i,j}$, let $M_2$ be obtained from $M_1$ by pulling out $\lambda^2$ as a common factor in each row of $M_1$ except the upper most one, and let $M_3$ be obtained from $M_2$ by pulling out $1/\lambda^2$ as a common factor in the left most column. Then $M_3=M$, and therefore
$$M_1 = \lambda^{2(n+1)} M_2 = \lambda^{2(n+1)-2} M_3 = \lambda^{2n} M.$$
Hence $M$ is either 0 or homogeneous of homogeneity degree $2n$. To see that $M \ne 0$, it is easy to see that, and it follows from Lemma 3.3 of \cite{impurity},
the determinant obtained from $M$ by plugging $x_{i,j}=1$ for every $x_{i,j}$ is $(-1)^{n-1} (n+1) \ne 0$. This completes the proof of (i).

(ii) The case $n=2$, i.e., the identity (\ref{Heron-M}), is easy to check. So we assume that $$n\ge 3.$$ 

Let $V_0$, $V_1$, and $V$ be the sets of indeterminates defined by
\begin{eqnarray}
V_1 = \{ x_{i,j} : 1 \le i < j \le n\},~V_0 = \{ x_{n+1,j} : 1 \le j \le n\},~
V = V_0 \cup V_1, \label{VV0V1}
\end{eqnarray}
and let $R_0$ and $R$ be the polynomial rings
\begin{eqnarray}
R_0 = k[V_0],~R = k[V]. \label{R0R}
\end{eqnarray}
For any non-zero $\phi \in R$,  let the total degree of  $\phi$ be denoted by $\deg (\phi)$, and let the degree of $\phi$ in the indeterminates $V_0$ be denoted by $\deg_0 (\phi)$. 
Since  the indeterminates $V_0$ appear only in the last row and last  column of $M$, it follows that
\begin{eqnarray}
\deg_0 (M) &\le& 4. \label{delta}
\end{eqnarray}

Now suppose  that $M$ is reducible, say 
\begin{eqnarray}
M &=& fg, \label{Mfg}
\end{eqnarray}
where $f$ and $g$ are homogeneous polynomials in $R$ of degrees $t, s \ge 1$.
We are to arrive at a contradiction.

For any $\phi \in R$, define $\phi^* \in R_0 [x]$  to be the polynomial 
obtained from $\phi$ by replacing every $x_{i,j}$ and  $x_{j,i}$, $1 \le i < j \le n$, by $x$. 
Then 
\begin{eqnarray}
M^* &=& f^*g^*. \label{Mfgstar}
\end{eqnarray}
By \cite[Theorem 4.1]{prekites}, $M^*$ is the Cayley-Menger determinant of the $n$-pre-kite $$PK[n;x;x_{n+1,1},\cdots,x_{n+1,n}]$$ and is given by
\begin{eqnarray} \label{Mstar=H}
M^*&=& (-x^2)^{n-2} H, \mbox{~where~} \nonumber\\
H &=& n(x^4+x_{n+1,1}^4+\cdots+x_{n+1,n}^4) -(x^2+x_{n+1,1}^2+\cdots+x_{n+1,n}^2)^2.
\end{eqnarray} 
By Theorem \ref{DDD}, $H$ is irreducible.
Since $M^* = f^* g^*$, it follows that $H$ divides one of the polynomials $f^*$ and $g^*$, say $g^*$. Therefore 
\begin{eqnarray}
\deg_0 (g) = \deg_0 (g^*) \ge \deg_0  (H) = 4,
\end{eqnarray}
and hence
\begin{eqnarray}
4 \ge \deg_0 (M) = \deg_0 (f) + \deg_0 (g) \ge  \deg_0 (f)   + 4.
\end{eqnarray}
Therefore $\deg_0 (f) = 0$, and hence $f$ does not contain any of the variables $x_{i,n+1}$, $1 \le i \le n$. Thus  $f$ must contain at least one of the other variables, say $x_{1,2}$. Let $F$ be obtained from $f$ by replacing $x_{1,2}$ by $x_{1,n+1}$. Since $M$ is symmetric, and since $f$ divides $M$, it follows that $F$ divides $M$, and therefore  $M=FG$, for some $G \in R$.
Thus $M^*=F^*G^*=(-x^2)^{n-2}H$. We now show that this is a contradiction.

Since $H$ is irreducible, it follows that either $H$ divides $F^*$ or $G^*$. The first case is impossible since  $F^*$ does not contain $x_{2,n+1}$. The second case implies that  $F^*$ divides $(-x^2)^{n-2}$, which is again impossible because $F^*$ contains  $x_{1,n+1}$. 

This completes the proof. \hfill $\Box$

\section{Remarks} \label{555} In the mathematical literature, there have been some polynomials which have attracted special attention, either for the elegance, or for their popping up in diverse, seemingly unrelated,  contexts. The polynomial $x^3+y^3+z^3-3xyz$ is the favorite polynomial alluded to in the title of \cite{favorite}, and it is the subject of Remark 4 (p.~193) of \cite{EM-ec}. Remark 5 (p.~194) of \cite{EM-ec} is devoted to the polynomial $x^3-(a^2-b^2-c^2)x+2abc$. 
 The Newton polynomial $x^3-(a^2+b^2+c^2)x+2abc$, too, has been a source of fascination, and has appeared in \cite{OMT}, and is the subject of \cite{NewtonDE}.
 The polynomials $f$ and $g$ defined in (\ref{f-a}) and (\ref{g}) above also  have their shares. The special case
\begin{eqnarray}
h=(a^2+x^2+y^2+z^2)^2-3(a^4+x^4+y^4+z^4)
\end{eqnarray}
which describes the relation among the distances of the vertices of an equilateral triangle of side length $a$ to an  arbitrary point in its plane is the tool that solves a problem that has appeared very frequently. The problem challenges readers to find $a$ given $x, y, z$ (or to find $z$ given $a$, $x$, $y$). The three known  numbers are usually given to be $3,4,5$, and the fourth unknown number is found to be some irrational  number that is not very pleasant, namely $\sqrt{25+12\sqrt{3}}$; see \cite[Problem 327, p. 162)]{500}, \cite[Chapter 1, Problem 27, p. 8; solution, p. 83--84]{Wagon}. A similar problem is solved in \cite[\S 55, p. 34]{Graham}, where the given numbers are $80$, $100$, and $150$, and where the answer is again not pleasant. However, the problem appears in \cite{CMJ}, where the given numbers are 5, 7, and $8$, and where the fourth number is found to be 3. One wonders whether there are integral  quadruplets other than $(3, 5, 7, 8)$. This naturally leads  to the Diophantine equation
\begin{eqnarray}\label{Sastry}
(w^2+x^2+y^2+z^2)^2-3(w^4+x^4+y^4+z^4)&=&0,
\end{eqnarray}
which is studied in \cite{Sastry},  where many (but not all) of its integer solutions are found. It would be an interesting 
challenge to try to find {\it all} the integer solutions of (\ref{Sastry}), and even of the more general (\ref{g}).


\end{document}